\documentclass[12pt]{article}
\usepackage{amssymb,latexsym,theorem,hyperref}

\newcommand{\Ab}{\mathop{\mathrm{Abutment}}\nolimits}
\newcommand{\Spec}{\mathop{\mathrm{Spec}}\nolimits}

\newcommand{\red}{{\mathop{\mathrm{red}}\nolimits}}
\newcommand{\ind}{\mathop{\mathrm{ind}}\nolimits}

\newcommand{\even}{{\mathrm{even}}}
\newcommand{\hgt}{\mathop{\mathrm{ht}}\nolimits}
\newcommand{\gr}{\mathop{\mathrm{gr}}\nolimits}
\newcommand{\im}{\mathop{\mathrm{im}}\nolimits}

\newcommand{\GL}{\mathop{\mathit{GL}}\nolimits}
\newcommand{\gl}{\mathop{\mathfrak{gl}}\nolimits}

\newcommand{\SL}{\mathop{\mathit{SL}}\nolimits}

\newcommand{\Z}{\mathbb Z}

\newcommand{\Gm}{{\mathbb G_m}}

\newcommand{\A}{{\cal A}}
\newcommand{\qed}{\unskip\nobreak\hfill\hbox{ $\Box$}}

\newtheorem{Theorem}{Theorem}[section]
\newtheorem{Proposition}[Theorem]{Proposition}
\newtheorem{Lemma}[Theorem]{Lemma}
\newtheorem{Corollary}[Theorem]{Corollary}
\theorembodyfont{\normalfont}
\newtheorem{Remark}[Theorem]{Remark}

\newtheorem{Definition}[Theorem]{Definition}

\begin{document}
\title{A reductive group with finitely generated cohomology algebras.}
\author{Wilberd van der Kallen}
\date{}
\maketitle
\sloppy
\begin{abstract}
Let $G$ be the linear algebraic group $\SL_3$ over a field $k$ of
characteristic two.
Let $A$ be a finitely generated commutative $k$-algebra on which $G$ acts
rationally by $k$-algebra automorphisms.
We show that the full cohomology ring $H^*(G,A)$ is finitely generated.
This extends the finite generation property of the ring of invariants
$A^G$. We discuss where the problem stands for other
geometrically reductive group schemes.
\end{abstract}

\section{Introduction}
Consider  a linear algebraic  group scheme
$G$ defined over a field $k$ of positive characteristic $p$.
So $G$ is an affine group scheme whose coordinate algebra $k[G]$ is finitely
generated as a $k$-algebra.
We say that $G$ has the cohomological finite generation property (CFG)
if the following holds.
Let $A$ be a finitely generated commutative $k$-algebra on which $G$ acts
 rationally by $k$-algebra automorphisms. (So $G$ acts on $\Spec(A)$.)
Then the cohomology ring $H^*(G,A)$ is finitely generated
as a $k$-algebra. Here, as in
\cite[I.4]{Jantzen}, we use the cohomology introduced by
Hochschild, also known as `rational cohomology'.

In \cite{cohGrosshans} we have shown that $\SL_2$ over a field of positive
characteristic has property (CFG). In this paper we explain how to modify the
argument so that it also covers $\SL_3$ in characteristic two.
We further explain what is missing to make further progress along the same
lines.

Recall that there is no analogous problem in
characteristic zero, because rational
cohomology vanishes in higher degrees for any
linearly reductive group.

\section{Geometric reductivity}
Let $G$ be a  linear algebraic
group scheme defined over the field $k$.
We first list a number of properties that $G$ may or may not have.
They are more familiar for linear algebraic groups than for linear algebraic
group schemes, but we will need them in the greater generality.

If $G$ acts rationally by $k$-algebra automorphisms on a commutative
$k$-algebra $A$, we will simply say that $G$ acts on $A$.

\paragraph{Property (FG)}
Whenever $G$ acts on a finitely generated $k$-algebra $A$, the
ring of invariants $A^G$ is a finitely generated $k$-algebra.

\paragraph{Property (Noeth)}Whenever $G$ acts on a finitely
generated $k$-algebra $A$, and $M$ is a noetherian $A$-module on which
$G$ acts compatibly (so the structure map $A\otimes M\to M$ is a map of
$G$-modules), the module of invariants $M^G$ is noetherian over $A^G$.

\paragraph{Property (Int)}Whenever $G$ acts on a finitely
generated $k$-algebra $A$, leaving invariant an ideal $J$, the ring
of invariants $(A/J)^G$ is integral over the image of $A^G$.

\begin{Remark}
In property (Int) one may drop the condition that $A$ is finitely generated.
The resulting
property is equivalent to the original.
\end{Remark}

\paragraph{Property (GR)}(Geometric reductivity)
Whenever $V$ is a finite dimensional $G$-module with one dimensional submodule
$L$ on which $G$ acts trivially, there
is a homogeneous $G$-invariant polynomial function
$f$ on $V$ with $f(0)=0$, $f|_{L}\neq0$.

\paragraph{Property (GRN)}Whenever $V$ is a finite dimensional $G$-module with
one dimensional submodule
$L$ on which $G$ acts trivially, the coordinate algebra $k[L]$ is a noetherian
$k[V]^G$-module. Here $k[V]$ is of course the coordinate algebra of $V$.

\paragraph{Property (GRI)}Whenever $V$ is a finite dimensional $G$-module with
one dimensional submodule
$L$ on which $G$ acts trivially, the coordinate algebra $k[L]$ is integral over
the image of
$k[V]^G$.

\paragraph{}We may now summarize some invariant theory as follows
\begin{Theorem}\label{GRequi}
The properties (FG), (Noeth), (Int), (GR), (GRN), (GRI) are equivalent.
\end{Theorem}
\paragraph{Proof}
That (GR), (GRN), (GRI) are equivalent is pretty obvious, as $k[L]$ is just
a polynomial ring in one variable, and the restriction map $k[V]\to k[L]$
is a surjective map of graded rings. It is a theorem of Nagata \cite{Nagata}
that
(GR) implies (FG). (By \cite{Borsari-Santos} the argument also works when
$G$ is a group scheme.) That (FG) implies   (Noeth) follows by considering
the symmetric algebra $S_A(M)$ on $M$ over $A$, or by considering its quotient
ring with underlying set $A\oplus M$. The rest is easier.
\qed

\paragraph{}
By $k$-algebra we always
mean an associative ring that contains $k$ in its center.
\begin{Lemma}
Let $A$ be a $k$-algebra and $l/k$ a field extension.
Then $A$ is finitely generated as a $k$-algebra if and only if
$A\otimes_k l$ is finitely generated as an $l$-algebra.
\end{Lemma}
\paragraph{Proof}Field extensions are faithfully flat.
\qed

\paragraph{}Encouraged by this lemma we further assume that $k$ is
algebraically
closed.

\begin{Lemma}
Let $G$ be a finite group scheme over $k$. Let it act
on a commutative $k$-algebra $A$. Then $A$ is integral over $A^G$.
In particular, $G$ is geometrically reductive.
\end{Lemma}
\paragraph{Proof}
In view of Theorem \ref{GRequi} the last statement may be left as an exercise.
When $G$ is reduced, hence just a finite group, the lemma is classical:
$a$ is a root of $\prod_{g\in G}(x-g(a))$. 
And if $G$ is not reduced, it has its first Frobenius kernel
$G_1$ as a normal subgroup and $A^G=(A^{G_1})^{G/G_1}$. Now the Lie algebra
\cite[I,7.7]{Jantzen}
acts by derivations on $A$
and $G_1$ acts trivially on the subalgebra $A^p$
generated by $p$-th powers. By induction on the rank $\dim_k(k[G])$
we may assume $A^{G_1}$ is integral over $A^G$.\qed

\subsection{Reductive}By Haboush \cite[II.10.7]{Jantzen} reductive groups are
geometrically reductive. By Popov \cite{PopovRed} a geometrically
reductive linear algebraic group has to be reductive, and by Waterhouse
\cite{Waterhouse}
a linear algebraic group scheme $G$ is geometrically
reductive exactly when the connected component $G_\red^o$ of
its reduced subgroup $G_\red$ is geometrically
reductive.

\section{Cohomological variations}

We keep our field $k$ algebraically closed of positive characteristic $p$.
We now introduce cohomological variants of the properties (FG),
(Noeth), (Int). We do not know how to make cohomological
variants of (GR), (GRI), (GRN).

\paragraph{Property (CFG)}
Whenever $G$ acts on a finitely generated $k$-algebra $A$, the cohomology
ring  $H^*(G,A)$ is a finitely generated $k$-algebra.

\paragraph{Property (CNoeth)}Whenever $G$ acts on a finitely
generated $k$-algebra $A$, and $M$ is a noetherian $A$-module on which
$G$ acts compatibly,
the cohomology $H^*(G,M)$ is noetherian as a module over $H^*(G,A)$.

\paragraph{Property (CInt)}Whenever $G$ acts on a finitely
generated $k$-algebra $A$, leaving invariant an ideal $J$, the even part
$H^\even(G,A/J)$ of the cohomology
ring  $H^*(G,A/J)$ is integral over the image of $H^\even(G,A)$.

\begin{Remark}
It is not essential to restrict to the even part, but the advantage of
that is that the even part is a commutative ring, so that we do not have
to explain the terminology `integral'. In any case, we will see that it is
the even part that truly matters. Note that $H^*(G,A)$ is a finitely
$k$-algebra if and only if  $H^\even(G,A)$ is a finitely generated $k$-algebra
and $H^*(G,A)$ is a noetherian module over it.
\end{Remark}

\begin{Remark}
In property (CInt) one may drop the condition that $A$ is finitely generated.
The resulting
property is equivalent to the original.
\end{Remark}

\begin{Lemma}
Property (CFG) implies properties (CNoeth) and (Cint).
\end{Lemma}
\paragraph{Proof}
Apply the same reasoning as in the proof of \ref{GRequi}.
\qed

\begin{Theorem}[Evens 1961 \cite{Evens}]\label{CFGEvens}
A finite group has property (CFG).
\end{Theorem}
Note that this is stronger than the theorem of Venkov \cite{Venkov},
 which refers only
to trivial coefficients. Our interest is in nontrivial
actions on finitely generated algebras.

After many attempts the result of Evens
was generalized by Friedlander and Suslin,
but they forgot to make it explicit in the same generality as Evens did.

\begin{Theorem}[Friedlander and Suslin 1997 \cite{Friedlander-Suslin}]
\label{CFGFS}
A finite group scheme has property (CFG).
\end{Theorem}

\paragraph{Proof} If $G$ is connected,
take $C=A^G$ in \cite[Theorem 1.5, 1.5.1]{Friedlander-Suslin}.
If $G$ is not connected, one  finishes the argument
by following \cite{Evens} as on pages 220--221 of
\cite{Friedlander-Suslin}. More specifically, on page 221 of
\cite{Friedlander-Suslin} one wants to replace `the subalgebra of $H^*(G,k)$
generated by the $\eta_i$' with `the $C$-subalgebra of  $H^*(G,C)$
generated by the $\eta_i$', where $C=A^G$ again. 
\qed

\begin{Lemma}
Let $G$ be a linear algebraic group and $H$ a
geometrically reductive subgroup scheme of $G$. Then $G/H$ is affine.
\end{Lemma}
\paragraph{Proof}
If $H$ is reduced, see Richardson \cite{Richardson}.
If not, choose $r$ so large that the image of $H$ under the $r$-th Frobenius
homomorphism $F^r:G\to G^{(r)}$ of \cite[I.9.4]{Jantzen} is reduced.
This image is then a geometrically reductive linear algebraic subgroup
of the linear algebraic group $G^{(r)}$.
The map $G/H\to F^r(G)/F^r(H)$ is finite \cite[III,\S
3,5.5b]{Demazure-Gabriel},
hence
affine, so $G/H$ is affine because
$F^r(G)/F^r(H)$ is affine.
\qed

\paragraph{}Now recall from \cite{cohGrosshans}
the following suggestive result.

\begin{Lemma}
Let $G$ be a linear algebraic group with property (CFG).
Then any geometrically reductive subgroup scheme $H$ of $G$ also
has property (CFG).
\end{Lemma}
\paragraph{Proof}
The  quotient $G/H$ is affine, hence there is no
cohomology of quasicoherent sheaves to worry about.
By \cite[I 4.6, I 5.13]{Jantzen} one simply gets,
when $H$ acts on a finitely generated $A$, that
$H^*(H,A)=H^*(G,\ind_H^{G}A)$.
And $\ind_H^{G}A=(k[G]\otimes A)^H$ is finitely generated.
\qed

\paragraph{} For instance, if $\GL_n$ as a group scheme over $k$ satisfies
(CFG), then so does $\SL_n$. Conversely, if $\SL_n$ satisfies (CFG),
then so does $\GL_n$, because
$H^*(\GL_n,A)=H^*(\SL_n,A)^\Gm$ for any $\GL_n$-algebra $A$.
Further every linear algebraic group scheme is a subgroup scheme of some
$\GL_n$---hence  of some $\SL_{n+1}$---and it is natural to ask

\subsection{Problem}
Prove that the linear algebraic group $\SL_{n}$ over $k$ has property (CFG).

\paragraph{} We now recall a principle that is a key ingredient
in the proofs of theorems \ref{CFGEvens} and \ref{CFGFS}. 
First a definition.
\begin{Definition}
We say that a spectral sequence $$E_r \Rightarrow \Ab$$ \emph{stops} if there is
an integer $r_0$ so that the differential $d_r:E_r\to E_r$ vanishes for $r\geq
r_0$.  
\end{Definition}

We now formulate the principle as a slogan.

\begin{Lemma}\label{noethss}
A noetherian spectral sequence stops.
\end{Lemma}

Let us explain what we mean by that.
Assume an associative ring $R$ acts on the levels of a spectral
sequence $$E_r \Rightarrow \Ab$$ in such a manner that each differential
$d_r:E_r\to E_r$ is an $R$-module map and each isomorphism
$\ker d_r/\im d_r \cong E_{r+1}$ is $R$-linear.
Say the spectral sequence
starts at level two.
Assume that $E_2$ is a noetherian $R$-module. Then the spectral sequence
stops. 
The proof is easy. One writes $E_r$ as $Z_r/B_r$ where $Z_r$, $B_r$ are the
appropriate submodules of $E_2$. As the $B_r$ form an ascending sequence,
we must have an $s$ so that $B_r=B_s$ for $r\geq s$. Then $d_r$ vanishes
for $r\geq s$. 

\paragraph{}
In the examples where the lemma is applied one often has a multiplicative
structure on the spectral sequence and the $R$-module structure of $E_2$
comes from a graded ring map $R\to E_2$.
Say $E_2$ is a finitely generated $k$-algebra, which we think of as given.
Then to show that $E_2$ is a noetherian $R$-module
one wants to see that the image of $R$ in $E_2$ is big enough. Indeed
the subring $E^\even_2\cap \mathrm{image(R)}$ should be so big that
$E^\even_2$ is integral over it.
Thus the
paradoxical situation is that in order to prove finite generation results (for
an abutment) one
needs to exhibit enough images (of elements of $R$ in $E_2$).
That is why we will be looking for universal cohomology classes (to take
cup product with).

\subsection{The Grosshans family}
{}From now on let $G$ be the linear algebraic group $\SL_n$ over $k$.
If $G$ acts on a finitely generated commutative $k$-algebra $A$,
Grosshans has studied in \cite{Grosshans contr}
a flat family over the affine line with general
fiber $A$ and special fiber $\gr A$. Here $\gr A$ is of course the associated
graded ring with respect to a certain filtration.
We write $\A$ for
 the finitely generated $k$-algebra whose spectrum
is the total space of the family. In characteristic zero the family had
been studied earlier
by Popov, and in characteristic $p$ a variant of $\gr A$
had been crucial in work
of Mathieu \cite{Mathieu G}. It was the more recent book
\cite{Grosshans book} that alerted us to possible relevance of
the family in the present
context.

\subsection{Good filtrations}
We choose a Borel group $B^+=TU^+$ of upper triangular matrices and the
opposite Borel group $B^-$. The roots of $B^+$ are positive.
If $\lambda\in X(T)$ is dominant, then $\ind_{B^-}^G(\lambda)$ is the
`dual Weyl module' or `costandard module'
$\nabla_G(\lambda)$ with highest weight $\lambda$. The formula
$\nabla_G(\lambda)=\ind_{B^-}^G(\lambda)$ just means that $\nabla_G(\lambda)$
is obtained from the Borel-Weil construction:
$\nabla_G(\lambda)$ equals $H^0(G/B^-,{\mathcal L})$ for a certain
 line bundle on the
flag variety $G/B^-$.
In a good filtration of a $G$-module the layers are traditionally
required to be of the form
$\nabla_G(\mu)$. However, to avoid irrelevant contortions when dealing
with infinite dimensional $G$-modules,
 it is important to allow a layer to be a direct sum
of any number of copies of the same $\nabla_G(\mu)$. We always follow that
convention (compare \cite[II.4.16 Remark 1]{Jantzen}).
We refer to \cite{vdkallen book} and
\cite{Jantzen} for proofs of the main properties of
good filtrations.
If $M$ is a $G$-module, and $m\geq-1$ is an integer so that
$H^{m+1}(G,\nabla_G(\mu)\otimes M)=0$ for all dominant $\mu$, then we say
as in \cite{Friedlander-Parshall}
that $M$ has \emph{good filtration dimension} at most $m$.
The case $m=0$ corresponds with $M$ having a good filtration.
And for $m\geq0$ it means that $M$ has a resolution
$$0\to M\to N_0 \to \cdots \to N_m\to 0$$ in which the $N_i$ have good
filtration.
We say that $M$ has good filtration dimension precisely $m$,
notation $\dim_\nabla(M)=m$,
 if $m$ is
minimal so that $M$ has good filtration dimension at most $m$.
In that case $H^{i+1}(G,\nabla_G(\mu)\otimes M)=0$ for all dominant $\mu$
and all $i\geq m$. In particular $H^{i+1}(G,M)=0$ for $i\geq m$.
If there is no finite $m$ so that  $\dim_\nabla(M)=m$, then we put
$\dim_\nabla(M)=\infty$.

\begin{Lemma}
Let $0\to M'\to M \to M''\to 0$ be exact.
\begin{enumerate}
\item $\dim_\nabla(M)\leq \max(\dim_\nabla(M'),\dim_\nabla(M''))$,
\item $\dim_\nabla(M')\leq \max(\dim_\nabla(M),\dim_\nabla(M'')+1)$,
\item $\dim_\nabla(M'')\leq \max(\dim_\nabla(M),\dim_\nabla(M')-1)$,
\item $\dim_\nabla(M'\otimes M'')\leq \dim_\nabla(M')+\dim_\nabla(M'')$.\qed
\end{enumerate}
\end{Lemma}

\subsection{Filtering $A$}
If $M$ is a $G$-module, and $\lambda$ is  a dominant weight,
then $M_{\leq\lambda}$ denotes the largest $G$-submodule all whose weights
$\mu$ satisfy $\mu\leq\lambda$ in the usual partial order
\cite[II 1.5]{Jantzen}. For example, $M_{\leq 0}=M^G$.
Similarly $M_{<\lambda}$ denotes the largest $G$-submodule all whose weights
$\mu$ satisfy $\mu<\lambda$. As in \cite{vdkallen book}, we form the
$X(T)$-graded module
$$\gr_{X(T)} M=\bigoplus_{\lambda\in X(T)}M_{\leq\lambda}/M_{<\lambda}.$$
We convert it to a $\Z$-graded
module through an additive height function $\hgt:X(T)\to \Z$ defined by
$\hgt=2\sum_{\alpha>0}\alpha^\vee$,
the sum being over the positive roots.
In other words, we put
$$M_{\leq i}=\sum_{\hgt(\lambda)\leq i}M_{\leq\lambda}$$
and then $\gr M$ is the associated graded module of the filtration
$M_{\leq 0}\subseteq M_{\leq 1}\cdots$.
We apply this in particular to our finitely generated commutative
$k$-algebra with $G$ action $A$.
This construction of $\gr A$ goes back to Mathieu \cite{Mathieu G}. From
Grosshans we learn to look at the specialization maps $\A\to \gr A$ and
$\A\to A$
(\cite[4.10, 4.11]{cohGrosshans}).

\subsection{Cohomology of the associated graded algebra}
Under the technical assumption, hopefully unnecessary, that either $n\leq5$
or $p>2^n$, we have shown in \cite[3.8]{cohGrosshans} that
there is a Hochschild-Serre spectral sequence
$$E_2^{ij}=H^i(G/G_r,H^j(G_r,\gr A))\Rightarrow H^{i+j}(G,\gr A)$$
which stops. Here $G_r$ is an appropriately chosen Frobenius kernel
in $G=\SL_n$.
(The choice of $r$ is not constructive. In contrast to
the work of Friedlander and Suslin we get only qualitative results.)
The fact that the spectral sequence stops was derived from an estimate of
good filtration dimensions, not by lemma \ref{noethss}.

{}From the stopping of the spectral sequence one gets

\begin{Theorem}\cite[Thm 1.1]{cohGrosshans}\label{fingrgr}
Let $\SL_n$ act on the finitely generated $k$-algebra $A$.
If $n<6$ or $p>2^n$, then $H^*(\SL_n,\gr A)$ is finitely generated as a
$k$-algebra.
\end{Theorem}

The filtration $A_{\leq 0}\subseteq A_{\leq 1}\cdots$ induces a filtration
of the Hochschild complex \cite[I.4.14]{Jantzen} whence a spectral sequence
$$E(A):E_1^{ij}=H^{i+j}(G,\gr_{-i}A)\Rightarrow H^{i+j}(G, A).$$
It lives in an unusual quadrant, but as long as the
$E_1$ is a finitely generated $k$-algebra this causes no difficulty with
convergence: given $m$ there will be only finitely many nonzero $E_1^{m-i,i}$. 
(Compare \cite[4.11]{cohGrosshans}. Note that in \cite{cohGrosshans}
the $E_1$ is mistaken for an $E_2$.)

Thus theorem \ref{fingrgr}
leads to the question whether the spectral sequence $E(A)$ also
stops. If it does, then we have no difficulty in transferring the finite
generation from $E_1$ to abutment. (As the differentials are derivations,
the $p$-th power of a homogeneous element of even degree
in $E_r$ automatically
passes to $E_{r+1}$. That makes that $E_{r+1}$ is finitely generated when
$E_r$ is.) As ring of operators to act on this spectral sequence
we may take $R=H^*(G,\A)$, with $\A$ the coordinate ring of the
Grosshans family (\cite[4.11, 4.12]{cohGrosshans}).
The trick is then to show that the natural map
$R\to E_1$ makes $E_1$ into a noetherian module.
If $\SL_n$ has property (CNoeth) or (CInt), then this is clear because $\gr A$
is a quotient of $\A$. We conclude

\begin{Theorem}
If $n<6$ or $p>2^n$, then properties (CFG), (CNoeth), (CInt) are equivalent
for $\SL_n$.
\end{Theorem}

To show that $E_1$ is a noetherian $R$-module
we wish to have some idea of the size of the image of $R$ in $E_1$,
and also of the size of $E_1$ itself.
Our $E_1$ is the abutment of the earlier Hochschild-Serre spectral sequence.
So we look back at that spectral sequence and try to get it also noetherian
over $R$, even though we already know it stops!
Looking at the proof of its stopping we see that actually the
Hochschild-Serre spectral sequence is noetherian over $H^\even(G,\gr A)$.
So it would follow from property (CInt) that the Hochschild-Serre
spectral sequence is indeed noetherian over $R$.
But let us look more directly at the image of $R^\even=H^\even(\SL_n,\A)$
in $H^0(G/G_r,H^*(G_r,\gr A))$.

\subsection{Two module structures}
Friedlander and Suslin did not just show that $H^*(G_r,\gr A)$ is a finitely
generated $k$-algebra, they actually provide an explicit finitely generated
graded $k$-algebra ${\cal S}^*$ with $G$ action
so that a natural graded map
$(\A)^{G_r}\otimes {\cal S}^*\to H^*(G_r,\gr A)$ makes the target into
a noetherian module over the source. (We suppress in the notation that
${\cal S}^*$ depends on $n$ and $r$.)
Here the subgroup scheme $G_r$ acts
trivially on ${\cal S}^*$, so that we also have an action of $G/G_r$ on
$(\A)^{G_r}\otimes {\cal S}^*$.
Then of course $H^0(G/G_r,H^*(G_r,\gr A))$ is also a noetherian
$H^0(G/G_r,(\A)^{G_r}\otimes {\cal S}^*)$-module.
Thus we like to compare the images of $R^\even$ and
$H^0(G,(\A)^{G_r}\otimes {\cal S}^*)=H^0(G/G_r,(\A)^{G_r}\otimes {\cal S}^*)$
in
$H^0(G/G_r,H^*(G_r,\gr A))$.
More specifically, we would like the image of $R^\even$ to contain the image
of $H^0(G,(\A)^{G_r}\otimes {\cal S}^*)$. That will prove that the
$E_2$ of the Hochschild-Serre spectral sequence
 is a noetherian $R$-module (cf. \cite[Prop. 3.8]{cohGrosshans}).
It suffices to factor the map $H^0(G,(\A)^{G_r}\otimes {\cal S}^{2m})\to
H^0(G/G_r,H^{2m}(G_r,\gr A))$ through $H^{2m}(G,\A)$. This is why
one would like to have universal classes in $H^{2m}(G,({\cal S}^{2m})^{\#})$,
where we use $V^\#$ as a notation for the linear dual of a vector space $V$.
Taking cup product with such a class yields a map from
$H^0(G,(\A)^{G_r}\otimes {\cal S}^{2m})$ to $H^{2m}(G,(\A)^{G_r})$ and
then one
has to show that this map fits in the obvious way in a commutative diagram.
This then necessitates putting
restrictions on the universal class.

\section{Universal cohomology classes}
Having explained why we care about having certain universal
cohomology classes, let us now get more specific. Further details are given in
\cite{cohGrosshans}.
Recall that $G=SL_n$ over an algebraically closed field $k$ of positive
characteristic $p$.

 Let $W_2(k)=W(k)/p^2W(k)$ be the ring of Witt vectors of length two over $k$,
see \cite[II \S6]{Serre}. If $V$ is a module for $G$ (or for $\GL_n$)
 we write $V^{(r)}$ for the
$r$-th Frobenius twist of $V$.
One has an extension of algebraic groups
$$1\to \gl_n^{(1)} \to \GL_n (W_2(k)) \to \GL_n(k)\to 1,$$
whence a cocycle class $e_1\in H^2(\GL_n,\gl_n^{(1)})$.
We call it the Witt vector class for $\GL_n$.
Its $m$-th cup power defines an element $e_1^{\cup m}$ in
$H^{2m}(G,(\gl_n^{(1)})^{\otimes m})$.
Now $(\gl_n^{(1)})^{\otimes m}$ contains the divided power
$\Gamma^m(\gl_n^{(1)})$ as the submodule of invariants under permutation
of the $m$ factors. We wish to lift $e_1^{\cup m}$ to a class in
$H^{2m}(G,\Gamma^m(\gl_n^{(1)}))$, in such a manner that some natural
properties are satisfied.
Recall that we have for $i\geq1$, $j\geq1$ an inclusion
$\Delta_{i,j}:\Gamma^{i+j}(\gl_n^{(1)})\to \Gamma^{i}(\gl_n^{(1)})
\otimes\Gamma^{j}(\gl_n^{(1)})$.

\begin{Definition}
We call the pair of integers $a$, $b$ \emph{special} if there is $i\geq0$
with $a=p^i$ and $1\leq b\leq (p-1)a$.
If $a$, $b$ is a special pair, then so is $ap^s$, $bp^s$ for any $s\geq0$.
Any integer larger than one can be written uniquely as $a+b$ with $a$,
$b$ special.
\end{Definition}

\subsection{Lifting Problem}\label{lift}
Put $c[1]=e_1$.
Show that there are
$c[m]\in H^{2m}(G,\Gamma^m(\gl_n^{(1)}))$ so that for every special pair
$a$, $b$ one has $$(\Delta_{a,b})_*(c[a+b])=c[a]\cup c[b]$$ in
$H^{2(a+b)}(G,\Gamma^a(\gl_n^{(1)})\otimes\Gamma^b(\gl_n^{(1)}))$.

\begin{Remark}\label{difference}
In \cite[Thm 4.4]{cohGrosshans} we required $(\Delta_{a,b})_*(c[a+b])=c[a]\cup
c[b]$
for all pairs $a$, $b$ with $a\geq1$, $b\geq1$.
\end{Remark}

\begin{Remark}\label{divres}
Look at the restriction of $c[m]$ to
$H^{2m}(G_1,\Gamma^m(\gl_n^{(1)}))=H^{2m}(G_1,k)\otimes \Gamma^m(\gl_n^{(1)})$.
The ring $\oplus_{m\geq0} H^{\even}(G_1,k)\otimes \Gamma^m(\gl_n^{(1)})$ is a
divided
power algebra over $H^\even(G_1,k)$ and
the restriction of $c[m]$ to $G_1$
is the $m$-th divided power of the restriction of $e_1$.
\end{Remark}

\paragraph{}
If one has solved the lifting problem and
$n\leq5$ or $p>2^n$, then one can construct the universal classes
in $H^{2m}(G,({\cal S}^{2m})^\#)$ alluded to
above and thus establish that $G$ has property (CFG).
This is what we did in \cite{cohGrosshans} for $G=\SL_2$, except that---as
mentioned in the remark \ref{difference}---we
gave the $c[m]$ some more properties than we require now.
So we must have another look at the proof of \cite[Lemma 4.7]{cohGrosshans}.
In this proof we must now take $a$, $b$ special.
The multiplication $S^a(\gl_n^{\#(r)})\otimes S^b(\gl_n^{\#(r)})\to
S^{a+b}(\gl_n^{\#(r)})$ is surjective, so that it still
follows by induction that the cup product with $c_{i}[a+b]^{(r-i)}$ takes the
desired values. As for the second part of the lemma, the restriction to
$ H^{2p^{r-1}}(G_{1},\gl_n^{(r)})$ is the same as before by remark
\ref{divres},
\cite[4.6]{cohGrosshans}. All in all the proofs in \cite{cohGrosshans}
generalize as soon as one has solved the lifting problem and
$n\leq5$ or $p>2^n$. But we do not have a method to decide the lifting problem
in general.

\paragraph{}
In our solution of the lifting problem for $\SL_2$ in
\cite{cohGrosshans} we used that the unipotent radical of a Borel
subgroup is abelian. That may be true for $\SL_2$,
but it fails for $\SL_3$, so we have to argue differently.
Put $\rho^\vee=\sum_{\alpha\mbox{ simple}}\alpha^\vee$, so that
$\rho^\vee(\varpi_i)=1$ for every fundamental weight $\varpi_i$.
The following proposition provides the missing link to deal with $\SL_3$
in characteristic two. It also gives a new proof for the $\SL_2$ case.

\begin{Proposition}\label{smalldim}
Let $n=2$ or let $n=3$ and $p=2$. Let $m\geq0$.
Let $M$ be a $G$-module with $\rho^\vee(\lambda)\leq m$ for every weight
$\lambda$ of $M$. Then the good filtration dimension $\dim_\nabla(M^{(1)})$
of $M^{(1)}$ is at most $m$.
\end{Proposition}

\begin{Corollary}
Let $n=2$ or let $n=3$ and $p=2$.
The maps $$(\Delta_{a,b})_*:H^{2(a+b)}(G,\Gamma^{a+b}(\gl_n^{(1)}))\to
H^{2(a+b)}(G,\Gamma^a(\gl_n^{(1)})\otimes\Gamma^b(\gl_n^{(1)}))$$ are
surjective
for $a\geq1$, $b\geq1$, so that the $c[m]$ of
problem \ref{lift} exist
and property (CFG) holds for $G$.
\end{Corollary}

\paragraph{Proof}The cokernel $C$ of $\Delta_{a,b}$
satisfies $H^{2(a+b)}(G,C)=0$, because
$\dim_\nabla(C)<2(a+b)$.\qed

\paragraph{Proof of Proposition}
This is rather straightforward. One starts with checking that the
cokernel of $\nabla(\varpi_i)^{(1)}\to \nabla(p\varpi_i)$ has a good
filtration for each fundamental
weight $\varpi_i$. This fails for other values of $n$, $p$.
Let us use induction on $m$.
The case $m=0$ of the proposition
is clear. Let $m>0$ and assume the result for strictly smaller
values. We may assume $M$ is finite dimensional and, taking a composition
series, we may assume $M$ is a simple module $L(\lambda)$ of high weight
$\lambda$ with $\rho^\vee(\lambda)= m$. Choose $\varpi_i$ so that
$\mu=\lambda-\varpi_i$ is dominant. Then
$\dim_\nabla(\nabla(\mu)^{(1)}\otimes \nabla(\varpi_i)^{(1)})\leq
\rho^\vee(\mu)+1\leq m$ and the cokernel of the embedding
$L(\lambda)^{(1)} \to \nabla(\mu)^{(1)}\otimes \nabla(\varpi_i)^{(1)}$
has good filtration dimension strictly less than $m$ by induction.\qed

\begin{Remark}
There is a reason why proposition \ref{smalldim} does not hold for other
combinations of $n$, $p$. We now sketch this reason.
If $V$ is the defining representation of $\GL_n$ then as another
corollary
one gets a surjection
$H^{4}(G,\Gamma^2(V^{(1)}))\otimes \Gamma^2(V^{(1)\#}))\to
H^{4}(G,(\gl_n^{(1)})^{\otimes 2})$.
In particular, let
$\tau\in H^{4}(G,\Gamma^2(V^{(1)}))\otimes \Gamma^2(V^{(1)\#}))$
lift $e_1\cup e_1$. Then cup product with $\tau$ describes a component
of the homomorphism $S^*((V\otimes V^\#)^\#)\to H^\even(G_1,k)$
that describes
the connection \cite[Thm 5.2]{Suslin-Friedlander-Bendel}
between cohomology of $G_1$ and the restricted nullcone.
(The restricted nullcone consists of the $n$ by $n$ matrices whose $p$-th
power is zero.)
But then the restricted nullcone can only contain matrices of rank at most one:
Two by two minors, viewed as elements of $S^2((V\otimes V^\#)^\#)$,
are annihilated by
$\Gamma^2(V))\otimes \Gamma^2(V^{\#})$, viewed as subspace of
$(S^2((V\otimes V^\#)^\#)^\#$.
That leaves only $n=2$ or $n=3$ and $p=2$.
\end{Remark}

\begin{Corollary}(cf. \cite[Cor. 4.13]{cohGrosshans})
Let $n=2$ or let $n=3$ and $p=2$.
Let $A$ be a finitely generated algebra with good filtration and let
$M$ be a noetherian $A$-module on which $G$ acts compatibly.
Then $M$ has finite good filtration dimension.
\end{Corollary}

\begin{Remark}
We find this surprising. We would have expected some smoothness condition
on $A$ or some flatness condition on $M$ to be needed. This makes one wonder
if a finitely generated algebra with good filtration is always a quotient
of a smooth finitely generated algebra with good filtration.
Note that the corollary forbids an algebra to be such a quotient when
the algebra does not have finite good filtration dimension.
In any case, this may be a good place to look when trying to defeat our
expectation that $\SL_n$ always satisfies (CFG). On the other hand, in a
later preprint (\texttt{arXiv:math.RT/0405238}) we show that the corollary extends
to the cases where $n<6 $ or $p>2^n$.
\end{Remark}

\paragraph{Proof of Corollary}We repeat the proof of
\cite[Cor. 4.13]{cohGrosshans}.
Let $k[G/U]$ denote the multicone $\ind_U^Gk$.
Recall that, for finite $d$,
 the good filtration dimension of $M$ is at most $d$ if and
only if $H^{d+1}(G,M\otimes k[G/U])$ vanishes. Now note that
$M\otimes k[G/U]$ is  a noetherian $A\otimes k[G/U]$-module and note that
$H^*(G,A\otimes k[G/U])$ is concentrated in degree zero.\qed

\end{document}